\documentclass[11pt,reqno]{amsart}
\usepackage{amsmath, amsthm, amssymb, amsfonts, mathrsfs, graphicx}
\usepackage{indentfirst}
\usepackage[mathscr]{euscript}
\usepackage{eqparbox}
\usepackage{appendix}

\newtheorem{theorem}{Theorem}[section]

\newtheorem{proposition}{Proposition}[section]
\numberwithin{equation}{section}

\def\pf{{\textit {Proof:} }}

\newcommand{\mysection}[1]{\section{#1}\setcounter{equation}{0}}

\newfont{\bb}{msbm10 at 11pt}

\newcommand{\bal}{\begin{aligned}}      \newcommand{\eal}{\end{aligned}}
\newcommand{\ba}{\begin{array}}      \newcommand{\ea}{\end{array}}
\newcommand{\bc}{\begin{center}}     \newcommand{\ec}{\end{center}}
\newcommand{\be}{\begin{enumerate}}  \newcommand{\ee}{\end{enumerate}}
\newcommand{\beq}{\begin{eqnarray}}  \newcommand{\eeq}{\end{eqnarray}}
\newcommand{\beQ}{\begin{eqnarray*}} \newcommand{\eeQ}{\end{eqnarray*}}
\newcommand{\bi}{\begin{itemize}}    \newcommand{\ei}{\end{itemize}}
\newcommand{\bt}{\begin{tabular}}    \newcommand{\et}{\end{tabular}}
\newcommand{\bdm}{\begin{displaymath}} \newcommand{\edm}{\end{displaymath}}

\newcommand{\Lrw}{\Longrightarrow}
\newcommand{\Llrw}{\Longleftrightarrow}

\def\qed{\hfill{Q.E.D.}\smallskip}

\newcommand{\ls}{\setlength{\baselineskip}{12pt}
	\setlength{\parskip}{3mm}}

\begin{document}
	
\title[The Dirac equation]{The Dirac equation on metrics of Eguchi-Hanson type}

\author[Z Cai]{Zhuohua Cai$^{\dag}$}
\address[]{$^{\dag}$School of Mathematics and Information Science, Guangxi University, Nanning, Guangxi 530004, PR China}
\email{2006301001@st.gxu.edu.cn}

\author[X Zhang]{Xiao Zhang$^{\flat}$}
\address[]{$^{\flat}$Guangxi Center for Mathematical Research, Guangxi University, Nanning, Guangxi 530004, PR China}
\address[]{$^{\flat}$Academy of Mathematics and Systems Science, Chinese Academy of Sciences, Beijing 100190, PR China}
\address[]{$^{\flat}$School of Mathematical Sciences, University of Chinese Academy of Sciences, Beijing 100049, PR China}
\email{xzhang@gxu.edu.cn, xzhang@amss.ac.cn}

\date{}

\begin{abstract}

We investigate parallel spinors on the Eguchi-Hanson metrics and find the space of complex parallel spinors are complex 2-dimensional. For the metrics of Eguchi-Hanson type with the zero scalar curvature, we separate variables for the harmonic spinors and obtain the solutions explicitly.\\\\
PASC numbers: 03.60.Pm, 04.20.Gz, 04.60.-m

\end{abstract}

\maketitle \pagenumbering{arabic}

\mysection{Introduction}\ls

Eguchi-Hanson metrics were constructed in \cite{EH1,EH2}, which are complete 4-dimensional Ricci flat, anti-self-dual ALE Riemannian metrics
on $R_{\geq 0} \times P^3$. They are one type of gravitational instantons, and play important roles in the Euclidean approach of gravitational quantization \cite{GH}. Since then, gravitational instantons including the Eguchi-Hanson metrics attract much attention both in geometry and in general relativity.

In order to investigate the physical behaviour of spin-$\frac{1}{2}$ particles in general relativity, Chandrasekhar solved the Dirac equation
in Kerr spacetimes by separating variables \cite{C}. We refer to \cite{FKSY, WZ} and references therein for further developments of the
the Dirac equation in various spacetimes. In \cite{SU}, Sucu and \"Unal solved the Dirac equation on Eguchi-Hanson and Bianchi $VII_0$ gravitational instanton metrics. In particular, the solutions are obtained as the product of two hypergeometric functions for the Eguchi-Hanson gravitational instanton metric.

As it actually relates to manifold's scalar curvature, it is interesting to study the Dirac equation on much more general metrics of Eguchi-Hanson type with the zero scalar curvature. This type of metrics were constructed by LeBrun using the method of algebraic geometry \cite{L} and by the second author solving an ordinary differential equation \cite{Z} on $R_{\geq 0} \times S^3 /Z_d$ $(d >2)$, and they provides the counter-examples of Hawking and Pope's generalized positive action conjecture.

In this note we first introduce the Eguchi-Hanson metrics, the metrics of Eguchi-Hanson type with the zero scalar curvature as well as the spin connections on them. Then we investigate the parallel spinors and the harmonic spinors. As the existence of a nontrivial parallel spinor on a Riemannian spin manifold implies that the manifold is Ricci-flat, parallel spinors only exist on the Eguchi-Hanson metrics, and we find the space of complex parallel spinors are complex 2-dimensional. For the metrics of Eguchi-Hanson type with the zero scalar curvature, we separate variables for the harmonic spinors and obtain the solutions explicitly.

\mysection{Metrics of Eguchi-Hanson Type}\ls

Metrics of Eguchi-Hanson type are given by
\begin{equation}
	g=f^{-2}\mathrm{d}r^{2}+r^{2}\left(\sigma_{1}^{2}+\sigma_{2}^{2}+f^2\sigma_{3}^{2}\right), \label{me2.1}
\end{equation}
where $ f $ is a real function of $r$, $\sigma_{1}$, $\sigma_{2}$ and $\sigma_{3}$ are the Cartan-Maurer forms for $ SU(2) \cong S^3 $
\begin{align*}
	&\sigma_{1}=\dfrac{1}{2}\left(\sin\psi \mathrm{d}\theta-\sin\theta \cos\psi \mathrm{d}\phi\right),\\
	&\sigma_{2}=\dfrac{1}{2}\left(-\cos\psi \mathrm{d}\theta-\sin\theta \sin\psi \mathrm{d}\phi\right),\\
	&\sigma_{3}=\dfrac{1}{2}\left(\mathrm{d}\psi+\cos\theta \mathrm{d}\phi\right)
\end{align*}
with three Euler angles $\theta$, $\phi$, $\psi$. The frame and coframe are
\begin{equation}
	e^{1}=f^{-1}\mathrm{d}r,\quad e^{2}=r\sigma_{1},\quad e^{3}=r\sigma_{2},\quad e^{4}=rf\sigma_{3}, \label{coframe}
\end{equation}
and
\begin{equation}
\begin{aligned}
		&e_{1}=f\partial_{r},\quad e_{2}=\dfrac{2}{r}\left(\sin\psi \partial _{\theta}-\dfrac{\cos\psi}{\sin\theta}\partial_{\phi}+\dfrac{\cos\theta \cos\psi}{\sin\theta}\partial_{\psi}\right),\\
		&e_{3}=\dfrac{2}{r}\left(-\cos\psi\partial_{\theta}-\dfrac{\sin\psi}{\sin\theta}\partial_{\phi}+\dfrac{\cos\theta \sin\psi}{\sin\theta}\partial_{\psi}\right),\quad e_{4}=\dfrac{2}{rf}\partial_{\psi}
	\end{aligned} \label{frame}
\end{equation}
respectively.

Let $\{ \omega^{i}_{j} \}$ be the connection 1-forms given by $\mathrm{d}e^{i}=-\omega^{i}_{j}\wedge e^{j}$. Then \cite{Z}
\begin{equation}
	\begin{aligned}
		&\omega^{2}_{1}=\dfrac{f}{r}e^{2},\quad \omega^{3}_{1}=\dfrac{f}{r}e^{3},\quad \omega^{4}_{1}=\left(\dfrac{f}{r}+f^{\prime}\right)e^{4},\\
		&\omega^{3}_{4}=\dfrac{f}{r}e^{2},\quad \omega^{4}_{2}=\dfrac{f}{r}e^{3},\quad \omega^{2}_{3}=\left(\dfrac{2}{rf}-\dfrac{f}{r}\right)e^{4}.
	\end{aligned} \label{connection}
\end{equation}

Let $B>0$ be certain positive constant. Taking
\begin{align*}
f=\sqrt{1-\dfrac{B}{r^4}},
\end{align*}
and
\begin{align*}
r \geq \sqrt[4]{B}, \quad 0 \leq \theta \leq \pi,\quad  0 \leq \phi \leq 2\pi,\quad  0 \leq \psi \leq 2\pi,
\end{align*}
metric (\ref{me2.1}) is the Eguchi-Hanson metric \cite{EH1, EH2}
\begin{equation}
	g=\left(1-\dfrac{B}{r^4}\right)^{-1}\mathrm{d}r^{2}+r^{2}\left(\sigma_{1}^{2}
	+\sigma_{2}^{2}+\left(1-\dfrac{B}{r^4}\right)\sigma_{3}^{2}\right), \label{me2.5}
\end{equation}
which is Ricci-flat and geodesically complete. Topologically, it is
\begin{align*}
R_{\geq 0} \times SU(2)/Z_2 \cong R_{\geq 0} \times SO(3) \cong R_{\geq 0} \times P_{3}.
\end{align*}
Let $d>2$ be a positive integer. Taking
\begin{align*}
f=\sqrt{1-\dfrac{2A}{r^2}-\dfrac{B}{r^4}}, \quad A=-\dfrac{d-2}{2}\sqrt{\dfrac{B}{d-1}}, \quad r_0=\sqrt[4]{\dfrac{B}{d-1}}
\end{align*}
and
\begin{align*}
r \geq r_0, \quad 0 \leq \theta \leq \pi,\quad  0 \leq \phi \leq 2\pi,\quad  0 \leq \psi \leq \frac{4\pi}{d},
\end{align*}
metric (\ref{me2.1}) is
\begin{equation}
	g=\left(1-\dfrac{2A}{r^2}-\dfrac{B}{r^4}\right)^{-1}\mathrm{d}r^{2}+r^{2}\left(\sigma_{1}^{2}+\sigma_{2}^{2}
	+\left(1-\dfrac{2A}{r^2}-\dfrac{B}{r^4}\right)\sigma_{3}^{2}\right), \label{me2.6}
\end{equation}
which is scalar-flat and geodesically complete, constructed in \cite{L, Z}. Topologically, it is
\begin{align*}
R_{\geq 0} \times SU(2)/Z_d \cong R_{\geq 0} \times S^3/Z_d.
\end{align*}

\mysection{Spin Connection and Parallel Spinors}\ls

We refer to \cite{H} for an introductory knowledge of spin geometry. It is clear that the underline manifolds of both metric (\ref{me2.5}) and metric (\ref{me2.6}) are spin and their complex spinor bundles are complex 4-dimensional equipped with spin connection
\begin{equation*}
	\nabla_{e_{k}}\Phi=e_{k}\Phi+\dfrac{1}{4}\sum_{i,j=1}^{4}g\left(\nabla_{e_{k}}e_{i},e_{j}\right)e_{i}\cdot e_{j}\cdot \Phi
\end{equation*}
for spinor $\Phi =\left(\Phi_{1}, \Phi_{2}, \Phi_{3}, \Phi_{4} \right)^t $. Using (\ref{connection}), we obtain
\begin{equation}
	\begin{aligned}
		&\nabla_{e_{1}}\Phi=e_{1}\Phi,\\
		&\nabla_{e_{2}}\Phi=e_{2}\Phi+\dfrac{f}{2r}\left(e_{1}\cdot e_{2}-e_{3}\cdot e_{4}\right)\cdot \Phi,\\
		&\nabla_{e_{3}}\Phi=e_{3}\Phi+\dfrac{f}{2r}\left(e_{1}\cdot e_{3}+e_{2}\cdot e_{4}\right)\cdot \Phi,\\
		&\nabla_{e_{4}}\Phi=e_{4}\Phi+\dfrac{1}{2}\left(\left(\dfrac{f}{r}+f^{\prime}\right)e_{1}\cdot e_{4}+\left(\dfrac{f}{r}-\dfrac{2}{rf}\right)e_{2}\cdot e_{3}\right)\cdot \Phi.
	\end{aligned} \label{spin d}
\end{equation}	

Throughout the paper, we fix the following representations for $e_i$ in $Spin(4)$.
\begin{equation}\label{repn}
\begin{aligned}
	&e_{1}\mapsto
	\begin{pmatrix}
		0&0&\eqmakebox[c]{1}&0\\
		0&0&0&\eqmakebox[c]{1}\\
		\eqmakebox[c]{-1}&0&0&0\\
		0&\eqmakebox[c]{-1}&0&0\\
	\end{pmatrix}
	,\qquad
	e_{2}\mapsto
	\begin{pmatrix}
		0&0&0&\eqmakebox[c]{i}\\
		0&0&\eqmakebox[c]{i}&0\\
		0&\eqmakebox[c]{i}&0&0\\
		\eqmakebox[c]{i}&0&0&0\\
	\end{pmatrix},\\
	&e_{3}\mapsto
	\begin{pmatrix}
		0&0&0&\eqmakebox[c]{-1}\\
		0&0&\eqmakebox[c]{1}&0\\
		0&\eqmakebox[c]{-1}&0&0\\
		\eqmakebox[c]{1}&0&0&0\\
	\end{pmatrix}
	,\qquad
	e_{4}\mapsto
	\begin{pmatrix}
		0&0&\eqmakebox[c]{i}&0\\
		0&0&0&\eqmakebox[c]{-i}\\
		\eqmakebox[c]{i}&0&0&0\\
		0&\eqmakebox[c]{-i}&0&0\\
	\end{pmatrix}.
\end{aligned}
\end{equation}

Spinor $ \Phi $ is parallel if $ \nabla_{X}\Phi=0 $ for any tangent vector $X$. Any Riemannian metric carried a parallel spinor must be Ricci-flat, see, e.g. \cite{H}. Thus it is interesting to investigate the parallel spinors for the Eguchi-Hanson metric.

\begin{theorem} \label{thm3.1}
Parallel spinors for the Eguchi-Hanson metric (\ref{me2.5}) are complex constant spinors which form a complex 2-dimensional space.
\end{theorem}
\pf Under (\ref{repn}), we have
	\begin{align}
		&\nabla_{e_{1}}\Phi=e_{1}
		\begin{pmatrix}
			\Phi_{1}\\
			\Phi_{2}\\
			\Phi_{3}\\
			\Phi_{4}\\
		\end{pmatrix}=0, \label{eqn3.3}\\
		&\nabla_{e_{2}}\Phi=e_{2}
		\begin{pmatrix}
			\Phi_{1}\\
			\Phi_{2}\\
			\Phi_{3}\\
			\Phi_{4}\\
		\end{pmatrix}
		+\dfrac{f}{2r}
		\begin{pmatrix}
			0\\
			0\\
			-2\mathrm{i}\Phi_{4}\\
			-2\mathrm{i}\Phi_{3}\\
		\end{pmatrix}
		=0, \label{eqn3.4}\\
		&\nabla_{e_{3}}\Phi=e_{3}
		\begin{pmatrix}
			\Phi_{1}\\
			\Phi_{2}\\
			\Phi_{3}\\
			\Phi_{4}\\
		\end{pmatrix}
		+\dfrac{f}{2r}
		\begin{pmatrix}
			0\\
			0\\
			2\Phi_{4}\\
			-2\Phi_{3}\\
		\end{pmatrix}
		=0, \label{eqn3.5}\\
		&\nabla_{e_{4}}\Phi=e_{4}
		\begin{pmatrix}
			\Phi_{1}\\
			\Phi_{2}\\
			\Phi_{3}\\
			\Phi_{4}\\
		\end{pmatrix}
		+\left(\dfrac{1}{rf}+\dfrac{f^{\prime}}{2}\right)
		\begin{pmatrix}
			0\\
			0\\
			-\mathrm{i}\Phi_{3}\\
			\mathrm{i}\Phi_{4}\\
		\end{pmatrix}=0. \label{eqn3.6}
	\end{align}
Firstly, (\ref{eqn3.3}) implies
\begin{align*}
\partial _r \Phi =0.
\end{align*}
Secondly, (\ref{eqn3.4}) gives
\begin{align*}
	\left(\sin\psi \partial _{\theta}-\dfrac{\cos\psi}{\sin\theta}\partial_{\phi}+\dfrac{\cos\theta \cos\psi}{\sin\theta}\partial_{\psi}\right)
	\begin{pmatrix}
		\Phi_{3}\\
		\Phi_{4}
	\end{pmatrix}
	=\dfrac{f}{r}
	\begin{pmatrix}
		\mathrm{i}\Phi_{4}\\
		\mathrm{i}\Phi_{3}
	\end{pmatrix}.
\end{align*}
Since the left hand side is independent on $r$, but the right hand side contains function $ f(r) $, it must be
\begin{align*}
\Phi_{3} =\Phi_{4} =0.
\end{align*}
Now (\ref{eqn3.6}) gives
\begin{align*}
\partial _\psi \Phi =0.
\end{align*}
Thus (\ref{eqn3.4}), (\ref{eqn3.5}) give
	\begin{align*}
		\dfrac{2}{r}\left(  \sin \psi \partial_{\theta}-\dfrac{\cos\psi}{\sin\theta}\partial_{\phi}\right)
		\begin{pmatrix}
			\Phi_{1}\\
			\Phi_{2}
		\end{pmatrix}
		 & =0,\\
\dfrac{2}{r}\left(  -\cos \psi \partial_{\theta}-\dfrac{\sin\psi}{\sin\theta}\partial_{\phi}\right)
		\begin{pmatrix}
			\Phi_{1}\\
			\Phi_{2}
		\end{pmatrix}
		 & =0.
	\end{align*}
Therefore
\begin{align*}
\partial_\theta \begin{pmatrix}
			\Phi_{1}\\
			\Phi_{2}
		\end{pmatrix}=
\partial_\phi \begin{pmatrix}
			\Phi_{1}\\
			\Phi_{2}
		\end{pmatrix}=0.
\end{align*}
So $ \Phi_{1} $ and $ \Phi_{2} $ are arbitrary complex constants and the proof of the theorem is complete. \qed

\mysection{Harmonic Spinors}\ls

Harmonic spinors are the solutions of the Dirac equation
\begin{align*}
\mathcal{D}=\sum_{k=1}  ^4 e_{k}\cdot \nabla_{e_{k}}=0.
\end{align*}
We shall find the explicit solutions of harmonic spinors on metric (\ref{me2.6}) of Eguchi-Hanson type with zero scalar curvature ($d>2$). In the case $d=2$, harmonic spinors on Eguchi-Hanson metric (\ref{me2.5}) were solved by Sucu and \"Unal in \cite{SU}, and the solutions can be obtained in terms of hypergeometric functions.

\begin{theorem}\label{Thm4.1}
Let positive integer $d >2$. There exists harmonic spinors
\beQ
\Phi=\mathrm{e}^{\mathrm{i}\left(n+\frac{1}{2}\right)\phi}
\begin{pmatrix}
	C_1 \mathrm{e}^{\mathrm{i}\frac{d}{2}\left(m+\frac{1}{2}\right)\psi}u_{1}(r)v_{1}(\theta)\\
	C_2 \mathrm{e}^{\mathrm{i}\frac{d}{2}\left(m-\frac{1}{2}\right)\psi}u_{2}(r)v_{2}(\theta)\\
	C_3 \mathrm{e}^{\mathrm{i}\frac{d}{2}\left(m+\frac{1}{2}\right)\psi}u_{3}(r)v_{3}(\theta)\\
	C_4 \mathrm{e}^{\mathrm{i}\frac{d}{2}\left(m-\frac{1}{2}\right)\psi}u_{4}(r)v_{4}(\theta)
\end{pmatrix},
\eeQ
on metric (\ref{me2.6}) of Eguchi-Hanson type with zero scalar curvature, where $m$, $n$ are integers, $C_i$ are complex constants,
\beq\label{uu}
\begin{matrix}
\left\{
\begin{aligned}
u_{1}(r)&=\Big(r^{2}-r_0 ^2 \Big)^{\frac{1-d-dm}{2d}}\Big(r^{2}+(d-1) r_0 ^2 \Big)^{-\frac{1-d-dm}{2d}-\frac{d\left(1+2m\right)}{4}},\\
u_{2}(r)&=\Big(r^{2}- r_0 ^2\Big)^{\frac{1-d+dm}{2d}}\Big(r^{2}+(d-1) r_0 ^2 \Big)^{-\frac{1-d+dm}{2d}-\frac{d\left(1-2m\right)}{4}},\\
u_{3}(r)&=\dfrac{1}{r}\Big(r^{2}-r_0 ^2\Big)^{-\frac{1}{2}-\frac{1-d-dm}{2d}}\Big(r^{2}+(d-1) r_0 ^2\Big)^{-\frac{1}{2}+\frac{1-d-dm}{2d}+\frac{d\left(1+2m\right)}{4}},\\
u_{4}(r)&=\dfrac{1}{r}\Big(r^{2}-r_0^2\Big)^{-\frac{1}{2}-\frac{1-d+dm}{2d}}\Big(r^{2}+(d-1) r_0^{2} \Big)^{-\frac{1}{2}+\frac{1-d+dm}{2d}+\frac{d\left(1-2m\right)}{4}},
\end{aligned}
\right.
\end{matrix}
\eeq
\beq\label{vv}
\begin{matrix}
\left\{
\begin{aligned} v_{1}(\theta)&=v_{3}(\theta)=\left(\sin\dfrac{\theta}{2}\right)^{n+\frac{1}{2}-\frac{d}{2}\left(m+\frac{1}{2}\right)}\left(\cos\dfrac{\theta}{2}\right)^{-n-\frac{1}{2}-\frac{d}{2}\left(m+\frac{1}{2}\right)},\\ v_{2}(\theta)&=v_{4}(\theta)=\left(\sin\dfrac{\theta}{2}\right)^{-n-\frac{1}{2}+\frac{d}{2}\left(m-\frac{1}{2}\right)}\left(\cos\dfrac{\theta}{2}\right)^{n+\frac{1}{2}+\frac{d}{2}\left(m-\frac{1}{2}\right)}.
\end{aligned}
\right.
\end{matrix}
\eeq
\end{theorem}
\pf Denote
\begin{align*}
F_{1}=\dfrac{f}{r}+\dfrac{f^{\prime}}{2}-\dfrac{1}{rf}, \quad F_{2}=\dfrac{2f}{r}+\dfrac{f^{\prime}}{2}+\dfrac{1}{rf}.
\end{align*}
Under (\ref{repn}), harmonic spinors satisfy
\begin{equation}\label{DE1}
	\begin{aligned}
&\left(f\partial_{r}-\dfrac{2\mathrm{i}}{rf}\partial_{\psi}+F_{1}\right)\Phi_{1}
-\dfrac{2}{r}\mathrm{e}^{\mathrm{i}\psi}\left(\partial_{\theta}-\mathrm{i}\csc\theta\partial_{\phi}+\mathrm{i}\cot\theta\partial_{\psi}\right)\Phi_{2}=0,\\	 &\left(f\partial_{r}+\dfrac{2\mathrm{i}}{rf}\partial_{\psi}+F_{1}\right)\Phi_{2}
+\dfrac{2}{r}\mathrm{e}^{-\mathrm{i}\psi}\left(\partial_{\theta}+\mathrm{i}\csc\theta\partial_{\phi}-\mathrm{i}\cot\theta\partial_{\psi}\right)\Phi_{1}=0,\\
&\left(f\partial_{r}+\dfrac{2\mathrm{i}}{rf}\partial_{\psi}+F_{2}\right)\Phi_{3}
+\dfrac{2}{r}\mathrm{e}^{\mathrm{i}\psi}\left(\partial_{\theta}-\mathrm{i}\csc\theta\partial_{\phi}+\mathrm{i}\cot\theta\partial_{\psi}\right)\Phi_{4}=0,\\
&\left(f\partial_{r}-\dfrac{2\mathrm{i}}{rf}\partial_{\psi}+F_{2}\right)\Phi_{4}
-\dfrac{2}{r}\mathrm{e}^{-\mathrm{i}\psi}\left(\partial_{\theta}+\mathrm{i}\csc\theta\partial_{\phi}-\mathrm{i}\cot\theta\partial_{\psi}\right)\Phi_{3}=0.
	\end{aligned}
\end{equation}

Separating variables
\begin{equation}\label{ST}
\begin{pmatrix}
	\Phi_{1}\\
	\Phi_{2}\\
	\Phi_{3}\\
	\Phi_{4}\\
\end{pmatrix}
=\mathrm{e}^{\mathrm{i}\left(n+\frac{1}{2}\right)\phi}
\begin{pmatrix}
	\mathrm{e}^{\mathrm{i}\frac{d}{2}\left(m+\frac{1}{2}\right)\psi}h_{1}(r,\theta)\\
	\mathrm{e}^{\mathrm{i}\frac{d}{2}\left(m-\frac{1}{2}\right)\psi}h_{2}(r,\theta)\\
	\mathrm{e}^{\mathrm{i}\frac{d}{2}\left(m+\frac{1}{2}\right)\psi}h_{3}(r,\theta)\\
	\mathrm{e}^{\mathrm{i}\frac{d}{2}\left(m-\frac{1}{2}\right)\psi}h_{4}(r,\theta)
\end{pmatrix},
\end{equation}
where $ m $, $ n $ are integers, the Dirac equation gives
\begin{equation}\label{SE}
	e^{\mathrm{i}\left(\frac{d}{2}-1\right)\psi}\vec{L}(r, \theta)=\vec{R}(r, \theta),
\end{equation}
where
\beQ
\begin{aligned}
\vec{L}(r, \theta) &=
\begin{pmatrix}
    \dfrac{rf}{2}\left(\partial_{r}+\dfrac{d\left(2m+1\right)}{2rf^{2}}+\dfrac{1}{r}+\dfrac{f^{\prime}}{2f}-\dfrac{1}{rf^{2}}\right)h_{1}\\
    \left(\partial_{\theta}-\left(n+\dfrac{1}{2}\right)\csc\theta+\dfrac{d}{2}\left(m+\dfrac{1}{2}\right)\cot\theta\right)h_{1}\\
    -\dfrac{rf}{2}\left(\partial_{r}-\dfrac{d\left(2m+1\right)}{2rf^{2}}+\dfrac{2}{r}+\dfrac{f^{\prime}}{2f}+\dfrac{1}{rf^{2}}\right)h_{3}\\
    \left(\partial_{\theta}-\left(n+\dfrac{1}{2}\right)\csc\theta+\dfrac{d}{2}\left(m+\dfrac{1}{2}\right)\cot\theta\right)h_{3}
\end{pmatrix},\\
\vec{R}(r, \theta) &=
\begin{pmatrix}
	\left(\partial_{\theta}+\left(n+\dfrac{1}{2}\right)\csc\theta-\dfrac{d}{2}\left(m-\dfrac{1}{2}\right)\cot\theta\right)h_{2}\\
	-\dfrac{rf}{2}\left(\partial_{r}-\dfrac{d\left(2m-1\right)}{2rf^{2}}+\dfrac{1}{r}+\dfrac{f^{\prime}}{2f}-\dfrac{1}{rf^{2}}\right)h_{2}\\
	\left(\partial_{\theta}+\left(n+\dfrac{1}{2}\right)\csc\theta-\dfrac{d}{2}\left(m-\dfrac{1}{2}\right)\cot\theta\right)h_{4}\\
	\dfrac{rf}{2}\left(\partial_{r}+\dfrac{d\left(2m-1\right)}{2rf^{2}}+\dfrac{2}{r}+\dfrac{f^{\prime}}{2f}+\dfrac{1}{rf^{2}}\right)h_{4}
\end{pmatrix}.
\end{aligned}
\eeQ

Note that $\vec{L}(r, \theta)$ and $\vec{R}(r, \theta)$ are independent with $\psi$, and $ d \geq 3 $, we obtain
\begin{align}
\vec{L}(r, \theta)=\vec{R}(r, \theta)=0. \label{LR0}
\end{align}
Thus, for $i=1, 3$, $j=2,4$, we have
\begin{align*}
\partial_{\theta}h_{i}&=\left(\left(n+\dfrac{1}{2}\right)\csc\theta-\dfrac{d}{2}\left(m+\dfrac{1}{2}\right)\cot\theta\right)h_{i}, \\ \partial_{\theta}h_{j}&=\left(-\left(n+\dfrac{1}{2}\right)\csc\theta+\dfrac{d}{2}\left(m-\dfrac{1}{2}\right)\cot\theta\right)h_{j}.
\end{align*}
Solving these ordinary differential equations, we obtain
\begin{align*} h_{i}&=u_{i}(r)\left(\sin\dfrac{\theta}{2}\right)^{n+\frac{1}{2}
-\frac{d}{2}\left(m+\frac{1}{2}\right)}\left(\cos\dfrac{\theta}{2}\right)^{-n-\frac{1}{2}-\frac{d}{2}\left(m+\frac{1}{2}\right)},\\
h_{j}&=u_{j}(r)\left(\sin\dfrac{\theta}{2}\right)^{-n-\frac{1}{2}
+\frac{d}{2}\left(m-\frac{1}{2}\right)}\left(\cos\dfrac{\theta}{2}\right)^{n+\frac{1}{2}+\frac{d}{2}\left(m-\frac{1}{2}\right)}.
\end{align*}
Substituting them into (\ref{LR0}), we have
\begin{align*}
	u_{1}^{\prime}&=\dfrac{1}{r^{4}-2Ar^{2}-B}\left(Ar-d\left(m+\dfrac{1}{2}\right)r^{3}\right)u_{1},\\
	u_{2}^{\prime}&=\dfrac{1}{r^{4}-2Ar^{2}-B}\left(Ar+d\left(m-\dfrac{1}{2}\right)r^{3}\right)u_{2},\\
	u_{3}^{\prime}&=\dfrac{1}{r\left(r^{4}-2Ar^{2}-B\right)}\left(3Ar^{2}+B+\left(d\left(m+\dfrac{1}{2}\right)-3\right)r^{4}\right)u_{3},\\
	u_{4}^{\prime}&=\dfrac{1}{r\left(r^{4}-2Ar^{2}-B\right)}\left(3Ar^{2}+B-\left(d\left(m-\dfrac{1}{2}\right)+3\right)r^{4}\right)u_{4}.
\end{align*}
Solving these ordinary differential equations, we obtain (\ref{uu}) and (\ref{vv}). This gives the proof of theorem. \qed

Finally we study the singular sets of harmonic spinors constructed in Theorem \ref{Thm4.1}. The singularity occurs when at least one of the power indices of $r-r_0$ is negative in (\ref{uu}) or at least one of the power indices of $\sin\dfrac{\theta}{2}$ or $\cos\dfrac{\theta}{2}$ is negative in (\ref{vv}). The relations between polar coordinates $\{r, \theta, \psi, \phi\}$ and Cartesian coordinates $\{x_0, x_1, x_2, x_3\}$ of the metric (\ref{me2.6}) are \cite{Z}
\begin{align*}
x_{1}&=r\cos\dfrac{\theta}{2}\cos\dfrac{\psi+\phi}{2}, \qquad x_{2}=r\cos\dfrac{\theta}{2}\sin\dfrac{\psi+\phi}{2},\\
x_{3}&=r\sin\dfrac{\theta}{2}\cos\dfrac{\psi-\phi}{2}, \qquad x_{0}=r\sin\dfrac{\theta}{2}\sin\dfrac{\psi-\phi}{2}.
\end{align*}
Since $r \geq r_0$, we obtain
\begin{align}
\sin \frac{\theta}{2}=0 \Llrw x_{3}=x_{0}=0, \quad \cos\frac{\theta}{2}=0 \Llrw x_{1}=x_{2}=0.\label{angu_sing}
\end{align}
Now we denote
\begin{align*}
a_m ^\pm & =\dfrac{1-d \pm dm}{2d}, \\
a_{mn} ^\pm & =n+\frac{1}{2}+\frac{d}{2}\Big(-m \pm \frac{1}{2}\Big), \\
b_{mn} ^\pm & =n+\frac{1}{2}+\frac{d}{2}\Big(m \pm\frac{1}{2}\Big).
\end{align*}
It is clearly
\begin{align}
a_{mn} ^+  -a_{mn} ^- =\frac{d}{2}, \quad b_{mn} ^+  -b_{mn} ^- =\frac{d}{2}. \label{abmn}
\end{align}

\begin{proposition}
The harmonic spinor $\Phi$, constructed in Theorem \ref{Thm4.1}, as well as its norm $|\Phi|$ are singular at $r=r_{0}$. Moreover, it holds that
\bi
\item[(1)] If $-\frac{d}{2}<a_{mn} ^- <0$, $-\frac{d}{2}<b_{mn} ^- <0$, $v_1$, $v_2$, $v_3$, $v_4$ and $|\Phi|$ are singular either at $x_1 x_2$-plane, or at $x_3 x_0$-plane;\\
\item[(2)] If $-\frac{d}{2}<a_{mn} ^- <0$, $b_{mn} ^- \geq 0$, either $v_1$, $v_2$, $v_3$, $v_4$, $|\Phi|$ are singular at $x_1 x_2$-plane, or $v_1$, $v_3$, $|\Phi|$ are singular at $x_3 x_0$-plane;\\
\item[(3)] If $-\frac{d}{2}<a_{mn} ^- <0$, $b_{mn} ^- \leq -\frac{d}{2}$, either $v_1$, $v_2$, $v_3$, $v_4$, $|\Phi|$ are singular at $x_1 x_2$-plane, or $v_2$, $v_4$, $|\Phi|$ are singular at $x_3 x_0$-plane;\\
\item[(4)] If $a_{mn} ^- \geq 0$, $-\frac{d}{2}<b_{mn} ^- <0$, either $v_2$, $v_4$, $ |\Phi| $ are singular at $x_1 x_2$-plane, or $v_1$, $v_2$, $v_3$, $v_4$, $|\Phi|$ are singular at $x_3 x_0$-plane;\\
\item[(5)] If $a_{mn} ^- \geq 0$, $b_{mn} ^- \geq 0$, $v_2$, $v_4$, $|\Phi|$ are singular at $x_1 x_2$-plane, or $v_1$, $v_3$, $|\Phi|$ are singular at $x_3 x_0$-plane,\\
\item[(6)] If $a_{mn} ^- \geq 0$, $b_{mn} ^- \leq -\frac{d}{2}$, $v_2$, $v_4$, $|\Phi|$ are singular either at $x_1 x_2$-plane or at $x_3 x_0$-plane;\\
\item[(7)] If $a_{mn} ^- \leq -\frac{d}{2}$, $-\frac{d}{2}<b_{mn} ^- <0$, either $v_1$, $v_3$, $|\Phi|$ are singular at $x_1 x_2$-plane, or $v_1$, $v_2$, $v_3$, $v_4$, $|\Phi |$ are singular at $x_3 x_0$-plane;\\
\item[(8)] If $a_{mn} ^- \leq -\frac{d}{2}$, $b_{mn} ^- \geq 0$, $v_1$, $v_3$, $|\Phi|$ are singular either at $x_1 x_2$-plane, or at $x_3 x_0$-plane;\\
\item[(9)] If $a_{mn} ^- \leq -\frac{d}{2}$, $b_{mn} ^- \leq -\frac{d}{2}$, either $v_1$, $v_3$, $|\Phi|$ are singular at $x_1 x_2$-plane, or $v_2$, $v_4$, $|\Phi|$ are singular at $x_3 x_0$-plane.
\ei
\end{proposition}
\pf The norm of harmonic spinors in Theorem \ref{Thm4.1} is
\begin{align*}
|\Phi|^{2}=\sum_{i=1}^{4}|\Phi_{i}|^{2}=\sum_{i=1}^{4}|C_{i}|^{2}u_{i}^{2}v_{i}^{2}.
\end{align*}
If any $u_{i}$ or $v_{i}$ goes to $\infty$ at some subset of the space, the harmonic spinors and their norm become singular there.

The power orders of $r-r_0$ in $u_i$ are
\begin{align*}
a_m ^-,  \quad a_m ^+, \quad -\frac{1}{2}-a_m ^-,  \quad -\frac{1}{2}-a_m ^+
\end{align*}
respectively. Note that $m$, $d$ are integers and $d>2$. Thus
\begin{align*}
m=0 & \Lrw a_m ^\pm <0, \quad -\frac{1}{2}-a_m ^\pm<0,\\
m \leq -1 & \Lrw a_m ^- >0, \quad  -\frac{1}{2}-a_m ^+ >0, \quad a_m ^+ <0, \quad -\frac{1}{2}-a_m ^-<0,\\
m \geq 1 & \Lrw a_m ^- <0, \quad -\frac{1}{2}-a_m ^+ <0, \quad a_m ^+ >0,\quad -\frac{1}{2}-a_m ^- >0.
\end{align*}
This indicates that $r_0$ is the singular point of the harmonic spinor $\Phi$ as well as its norm $|\Phi|$.

From (\ref{vv}), we obtain
\begin{align*}
v_{1}(\theta)&=v_{3}(\theta)=\left(\sin\frac{\theta}{2}\right)^{a^- _{mn}}
                             \left(\cos\frac{\theta}{2}\right)^{- b^- _{mn} -\frac{d}{2}},\\
v_{2}(\theta)&=v_{4}(\theta)=\left(\sin\frac{\theta}{2}\right)^{-a^- _{mn} -\frac{d}{2}}
                             \left(\cos\frac{\theta}{2}\right)^{b^- _{mn}}.
\end{align*}
Therefore the rest parts of the theorem follow from the above formulas, (\ref{angu_sing}) and the fact that $\sin\frac{\theta}{2}$, $\cos\frac{\theta}{2}$ can not be both zero. \qed

\bigskip

{\footnotesize {\it Acknowledgement. The authors would like to thank referee's valuable suggestions. This work is supported by Chinese NSF grants 11731001, the special foundations for Guangxi Ba Gui Scholars and Junwu Scholars of Guangxi University.}


\begin{thebibliography}{99}

\bibitem{EH1} T. Eguchi, A.J. Hanson, Asymptotically flat self-dual solutions to Euclidean gravity, Phys. Lett. \textbf{74B} (1978) 249-251.
\bibitem{EH2} T. Eguchi, A.J. Hanson, Self-dual solutions to Euclidean gravity,  Ann. Phys. \textbf{120} (1979) 82-106.
\bibitem{GH} G.W. Gibbons, S.W. Hawking, Classification of gravitational instanton symmetries, Commun. Math. Phys. \textbf{66} (1979) 291-310.
\bibitem{C} S. Chandrasekhar, The solution of Dirac's equation in Kerr geometry, Proc. R. Soc. A \textbf{349} (1976) 571-575.
\bibitem{FKSY} F. Finster, N. Kamran, J. Smoller, S.T. Yau, Nonexistence of time-eriodic solutions of the Dirac equation in an axisymmetric black hole geometry, Commun. Pure Appl. Math. \textbf{53} (2000) 902-929.
\bibitem{WZ} Y.H. Wang, X. Zhang, Nonexistence of time-periodic solutions of the Dirac equation in non-extreme Kerr-Newman-AdS spacetime, Sci. China Math. \textbf{61} (2018) 73-82.
\bibitem{SU} Y. Sucu, N. \"Unal, Dirac equation in Euclidean Newman-Penrose formalism with applications to instanton metrics, Class. Quantum Grav. \textbf{21} (2004) 1443-1451.
\bibitem{L} C. LeBrun, Counter-examples to the generalized positive action conjecture, Commun. Math. Phys. \textbf{118} (1988) 591-596.
\bibitem{Z} X. Zhang, Scalar flat metrics of Eguchi-Hanson type, Commun. Theor. Phys. (Beijing, China) \textbf{42} (2004) 235-238.
\bibitem{H} O. Hijazi, Spectral properties of the Dirac operator and geometrical structures, Proceedings of the Summer School on Geometric Methods in Quantum Field Theory, Villa de Leyva, Colombia, July 12-30, 1999, World Scientific, Singapore (2001), 116-169.
	
\end{thebibliography}
\end{document}